# Performance Comparison of PSO Based State Feedback Gain (K) Controller with LQR-PI and Integral Controller for Automatic Frequency Regulation

NARESH KUMARI[1], A. N. JHA[2], NITIN MALIK[3]

[1,3]School of Engineering and Technology
The NorthCap University,
Sector 23A, Gurgaon, Haryana, India
nareshkumari@ncuindia.edu, nitinmalik77@gmail.com

[2] Ex-Professor, Electrical Engineering Department, Indian Institute of Technology
New Delhi, India
anjha@ee.iitd.ac.in

*Abstract:* In the present work a new controller called Particle Swarm Optimization (PSO) based state feedback gain (K) controller has been proposed for frequency regulation of a two area system and then its performance is compared with earlier designed controllers such as Linear Quadratic Regulator–Proportional Integral (LQR-PI) controller and Integral controller. The performance comparison has been done for the power system network comprising of two thermal power plants which are tie line connected. For using the optimal control based method such as LQR-PI controller and computationally intelligent method such as PSO based state feedback gain (K) controller, the state space modeling of the system has been done. Transfer function model for the system is used for finding the response of Integral controller. In an effective generation control scheme the change in frequency should be minimum during the load variation. The proposed PSO based state feedback gain (K) controller technique has been found most effective for improving the frequency response.

## 1. Introduction

The electrical power system has the main task of keeping a match between the power generation with the load requirement which is necessary to supply the users with affluent, quality and stable power. The power system has to maintain the optimum values of frequency and voltage level during the different load perturbations. With the change in the power demand there is a tendency of change in the frequency of the interconnected power system along with the power exchange between different areas. Further the load frequency control (LFC) can be achieved by proper regulation of the participating generating units. The area control error (ACE) has to be minimized for the frequency enhancement and the integral square error is generally taken as ACE. The objective of automatic load frequency controller (ALFC) is to reduce the integral square error (ISE) to zero with the continuous change of demand of active power so that the total generated power of the system and load requirement properly match with each other [1]. The control of active power for LFC is possible by varying the parameters of various controllers and sources in the system. There are basically three types of possible connections of controllers with the power system network such as central controller, single agent and decentralized controller [2]. In case of central controller only one control unit is used and it is connected to all the areas in the power system network to control their parameters and maintain the generation and demand balance. Further in single agent method the control action is performed on individual distributed generation system. In case of decentralized control each area is equipped with individual controller and local feedback control is used. The proposed PSO based state feedback gain (K) controller is used in central control mode and one controller is used for both the areas of power system network.

In this work the thermal network with two-area systems has been designed using the reheat turbines by developing the transfer function model in MATLAB/SIMULINK.. The total number of states taken to develop the model of the system are eleven [3]. Three load frequency control techniques such PSO based state feedback gain (K) controller, LQR-PI controller and Integral controller have been applied each time on the developed model. Then the performance comparison is done for three types of load frequency controllers. The PI controller can be tuned by using the optimal behavior of Linear Quadratic Regulator in case of LQR-PI controller. The set point tracking is possible with the combination of state feedback based controller and proportional integral controller as in such type of controller the advantages of both types of controllers can be attained [4]. The power system transfer function model has been developed by combining the transfer function for different parts of thermal power network such as generator, governor, turbine, reheat and tie line power. Further, the load frequency enhancement has also been done through optimization of state feedback gain (K) with PSO.

## 2. Power System Modelling

The system investigated is designed for two interconnected thermal power plants. Two power system areas are connected to each other by a tie-line and the reheat turbine is used in both the power plants. Each thermal plant taken for performance comparison of different load frequency controllers is of generation capacity 2000 MW. The system model developed is made linear by making certain assumptions and approximations in the mathematical model. The model based on transfer function approach has been designed in MATLAB environment as well as state space modeling of system has been done to analyze the performance of Integral action based controller and LQR-PI controller [5].

The power system state variables have been determined to obtain the set of equations of the interconnected thermal power network in state space. Two area based thermal system block diagram with eleven state variables for automatic frequency regulation of a linear network is shown Fig. 1.The model has been developed taking reheat into consideration for both the thermal power systems.

All state feedback variables of the system as shown in Figure 2 are as below:

$x1 = \Delta f1$, $\quad x2 = \Delta Pt1$,
$x3 = \Delta Pr1$, $\quad x4 = \Delta Pg1$
$x5 = \Delta f2$, $\quad x6 = \Delta Pt2$,
$x7 = \Delta Pr2$, $\quad x8 = \Delta Pg2$
$x9 = \Delta Pt\,tie(1,2)$, $x10 = \int (ACE1)dt$, $x11 = \int (ACE2)dt$

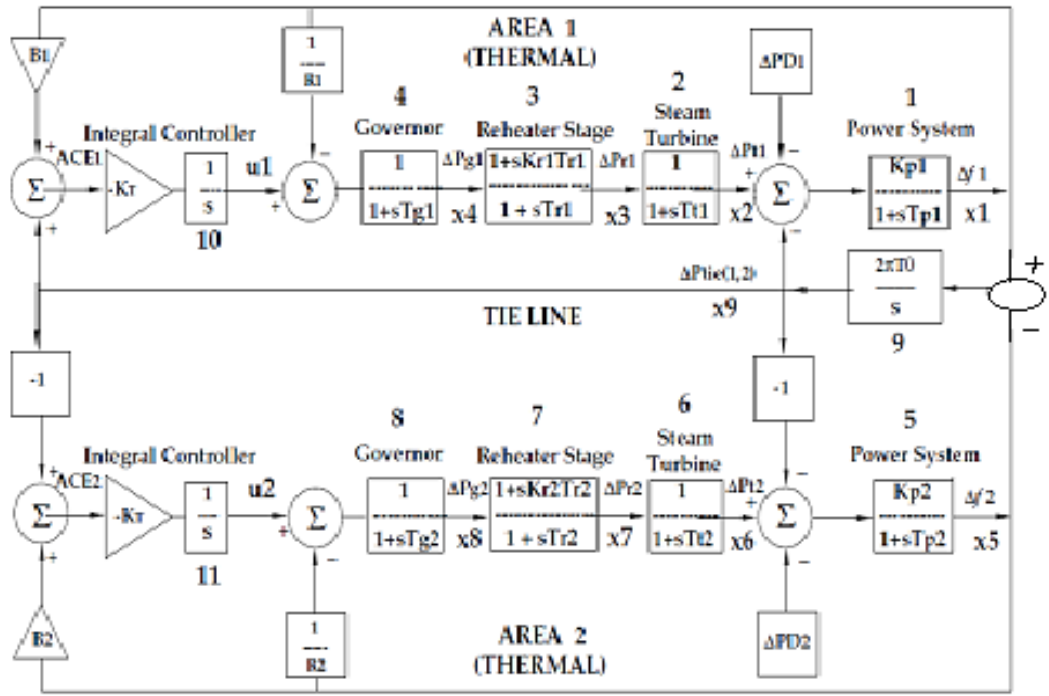

Fig. 1 Two area thermal power network represented with transfer function

The state vector X is comprising of all eleven states and control vector with control inputs u1 and u2.

The model of the system is designed in MATLAB environment on the basis of state space considering all states of two reheat thermal power systems as shown in Fig. 2 for analysis of the performance of LQR-PI controller. The PI controller performance is analyzed with transfer function model. Further for the analysis of PSO based feedback gain (K) controller the state space model is used.

The equation for optimal control of the system power is:

$u(t) = -Kx(t)$ ---------------------------- (1)

where state feedback gain (K) can be given as :

$K = R^{-1} B^T P$ -------------------------------- (2)

The matrix P can be find out by solving Riccati equation given by (3):

$A^T P + PA - P BR^{-1} B^T P + Q = 0$ ---- (3)

The system Performance Index in case of LQR-PI controller is as follows:

Performance Index = $0.5 \int ( x'Q x + u'R u) dt$ - ---------------------------------------------------- (4)

Further the performance index for PSO based state feedback gain (K) controller can be given as an integral square of error (ISE) as given below:

$ISE = 1/2 \int (\Delta f_1^2 + \Delta f_2^2) dt$ ------------------(5)

Both the matrices related to the weights such as Q and R are chosen judiciously for the stability of the system and for proper designing of LQ regulator. The optimal regulators are designed by changing the diagonal elements of Q and the elements of matrix R are not changed [6]. The weight matrices are also selected according to the predefined specifications of the system and desired response of the system. For PSO based feedback gain (K) controller the initial values of Q and R are selected as unity for both the diagonal matrices in order to reduce the complexity of the system.

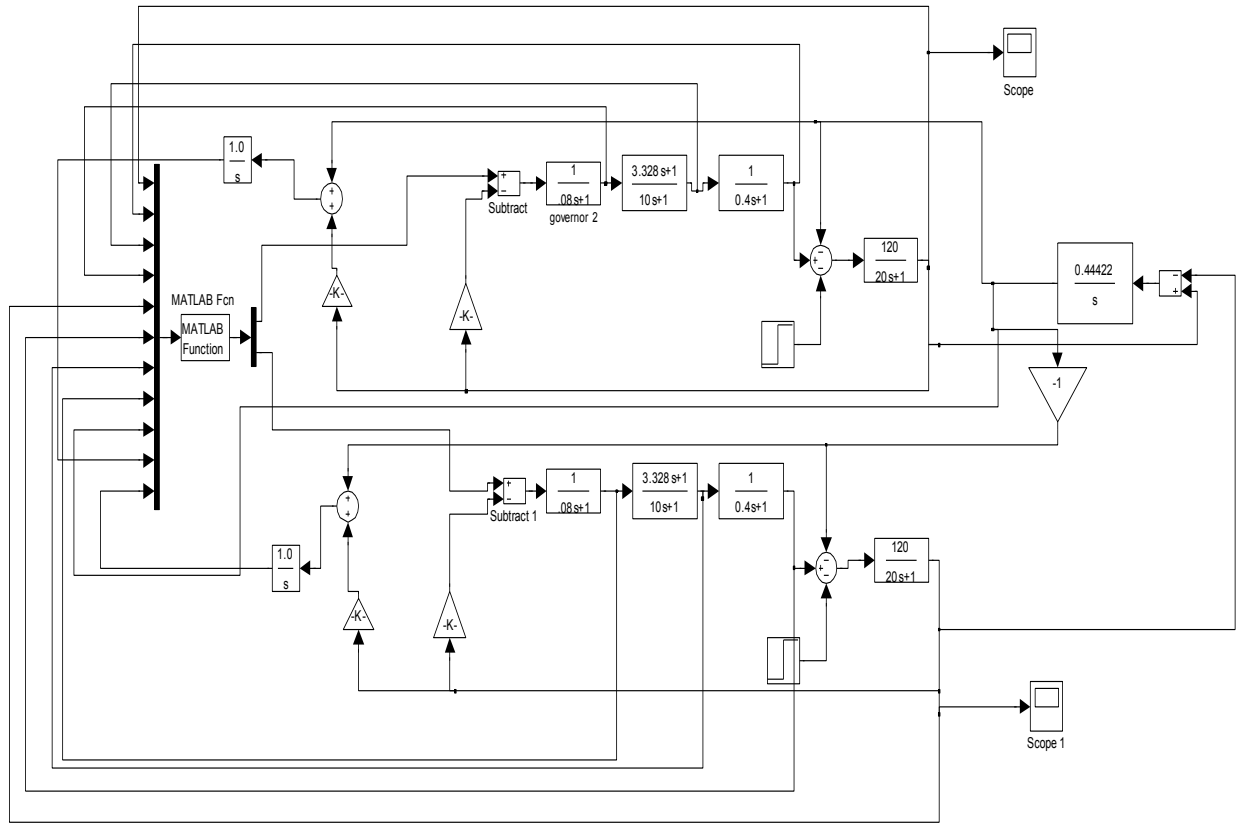

Fig. 2 Two area interconnected power system network along with reheat turbine for using LQR based PI control simulated in MATLAB/SIMULINK [5]

## 3. Controller Design for ALFC

The performance of the three types of controllers has been compared in this work for the automatic frequency regulation of the interconnected two equal thermal power systems with each of the rating of 2000 MW. The system response has been observed by changing the loading in area 1 by 1% and it has been observed that this will affect the frequency response of both the areas. In the same manner the loading has been changed in area 2 and its effect has been observed on the system. In this section, a description is given for three types of controllers viz. Integral controller, LQR- PI controller and PSO based state feedback gain (K) controller for automatic frequency regulation of interconnected power generation network. The performance comparison of these controllers is described in section 4 with their frequency response plot.

**3.1** The first controller used for frequency response enhancement is an Integral controller. The integral gain (Ki) is tuned by PSO technique to achieve the robust performance. As there are certain structural limits for the parameters which have to be considered for designing of a controller, the probabilistic search technique such as particle swarm optimization (PSO) is used to find the control parameter of the integral load frequency controller within the structural constraints [7].

**3.2** In the second method for the AGC the Linear Quadratic Regulator is used to optimize the integral controller gain as well as the states of system. First algebraic Riccati equation has been used to solve for the matrix P with equation (3). The performance index in this method is given in equation (4). For designing of weight matrices Q and R the proper understanding of the system is essential. Q and R are mainly the diagonal matrices. The matrices A, B and C for the system are written by observing the state differential equations from equations and specifications of the system. These matrices for the rating of 2000 MW two power plants can be given as:

$$A = \begin{bmatrix} -0.05 & 6 & 0 & 0 & 0 & 0 & 0 & 0 & -6 & 0 & 0 \\ 0 & -3.33 & 3.33 & 0 & 0 & 0 & 0 & 0 & 0 & 0 & 0 \\ -2.08 & 0 & -0.1 & -6.1 & 0 & 0 & 0 & 0 & 0 & 0 & 0 \\ -4.16 & 0 & 0 & -12.5 & 0 & 0 & 0 & 0 & 0 & 0 & 0 \\ 0 & 0 & 0 & 0 & -0.05 & 6 & 0 & 0 & 6 & 0 & 0 \\ 0 & 0 & 0 & 0 & 0 & -3.33 & 3.33 & 0 & 0 & 0 & 0 \\ 0 & 0 & 0 & 0 & -2.08 & 0 & -0.1 & -6.1 & 0 & 0 & 0 \\ 0 & 0 & 0 & 0 & -4.16 & 0 & 0 & -12.5 & 0 & 0 & 0 \\ 3.42 & 0 & 0 & 0 & -3.42 & 0 & 0 & 0 & 0 & 0 & 0 \\ 2 & 0 & 0 & 0 & 0 & 0 & 0 & 0 & 1 & 0 & 0 \\ 0 & 0 & 0 & 0 & 2 & 0 & 0 & 0 & -1 & 0 & 0 \end{bmatrix}$$

$$B = \begin{bmatrix} 0 & 0 \\ 0 & 0 \\ 6.250 & 0 \\ 12.50 & 0 \\ 0 & 0 \\ 0 & 0 \\ 0 & 6.25 \\ 0 & 12.5 \\ 0 & 0 \\ 0 & 0 \\ 0 & 0 \end{bmatrix}$$

C = [1 0 0 0 1 0 0 0 0 0 0];

In the MATLAB environment the state space model is given by the equation (6)

$$sys = ss(A,B,C,D) \text{ -------------------- (6)}$$

The matrix K can be calculated with LQR method as follows:

$$[K] = lqr(sys, Q, R) \text{ -------------------- (7)}$$

Further the value of matrix K obtained with LQR-PI controller for the system under study is given as:

K = [0.4896 1.0873 2.3995 0.7172 0.026
0.0026 0.1621 0.053 0.823 1 0 ;
0.026 0.0026 0.1621 0.053 0.4896
1.0873 2.3995 0.7172 0.823 0 1] ------(8)

**3.3** In the PSO based state feedback gain (K) method the elements of matrix K act as the particles [8, 9]. There are many other computationally intelligent techniques for optimizing the controller parameters as in literature [10, 11, 12]. The system investigated in this work is linearized, some of the previous work in LFC has taken the non-linear models of the power system network [13, 14, 15].

The pseudo code using PSO for the calculation of optimized value of state feedback gain matrix (k) is given below:

For every particle:

Initialize

Do

Solve the equation (2) for obtaining the initial values of matrix (K)

Calculate the particle velocity given by

$$V_i(t) = V_i(t-1) + c1\,\varphi 1(p\,(p\,best - x_i\,(t-1)) + c2\,\varphi 2\,(p\,(gbest - x_i\,(t-1)) \text{ ------------- (9)}$$

Update particle position according to equation (10):

$$x_i(t) = x_i(t-1) + v_i(t) \text{ ----------- (10)}$$

Till the stopping criterion is met.

End

The present value of ACE is compared to best value of ACE in previous iterations (pbest).

Then this new pbest is considered as the current value for the pbest.

The particle at which the least ACE is obtained taken as gbest.

The optimum values of elements of matrix K of dimension 2X11 are obtained with PSO technique. The model on which the PSO technique has been applied is described in section II. This is a fast and efficient tuning approach. The population size taken in PSO method is 100, search range is taken as [-3, 3], acceleration constants are chosen as c1=1.5, c2=1.5 and iteration number is taken as 100. As the main advantage of PSO is that the initial value of the parameter to be optimized can be taken as zero and the suitable range for optimum value is given for fast convergence. So initially each element of K matrix has been taken as zero.

As per the system and its loading conditions the final K matrix after applying PSO method is:

K = K'T ---------------------------------(11)
K = [0.012   0.0707   2.8555   -0.7412   0.7286   0.6396   0.4685   -0.1815   1.3475 1.0000 0.0000;
    0.7286   0.6396   0.4685   -0.1815   1.1724   2.0707   2.8555   0.7412   1.3475 0.0000 1.0000]------------------------(12)

Finally the control vector can be obtained as per the equation (1).The performance index of the system for finding the optimum state feedback gain (K) matrix using PSO technique is given in equation (5).

## 4. Simulation Results and Performance Comparison of Proposed PSO based K Controller with LQR-PI and Integral Controllers

The model for which the performance of the three controllers has been compared is explained in section II. Three types of controllers have been used for frequency response enhancement of the system considering one controller at a time. The various parameters of the controllers have been optimized as described in earlier section.

The loading has been increased by 1% at a time in one thermal power area for the simulation study done in MATLAB environment. It has been observed that the load change in area-I results in the frequency variation in the area-I as well as in area-II. Then the performance of controllers is analyzed by increasing the loading by same amount in area II. The optimally designed PI controller, LQR-PI controller and PSO based state feedback gain (K) matrix controller helps in keeping all the state variables in proper range along with the permissible frequency deviation during load change.

The change in frequency for both the areas for loading in area I have been shown in Fig. 3-Fig. 4.The frequency responses are compared for three types of controllers such as Integral, LQR-PI controller and PSO based K .The loading effect in area II and thus frequency responses with all controllers have been compared in Fig. 5-Fig.6.

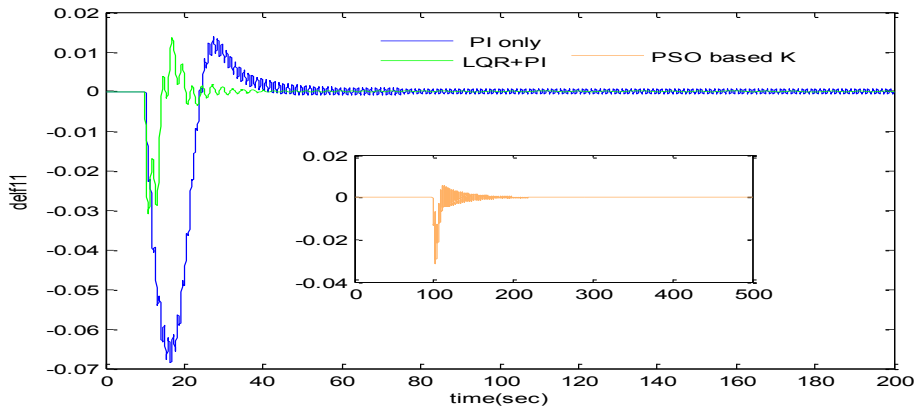

Fig. 3: Δf1  Vs t  for area 1 if 1% load variation occurs in area 1

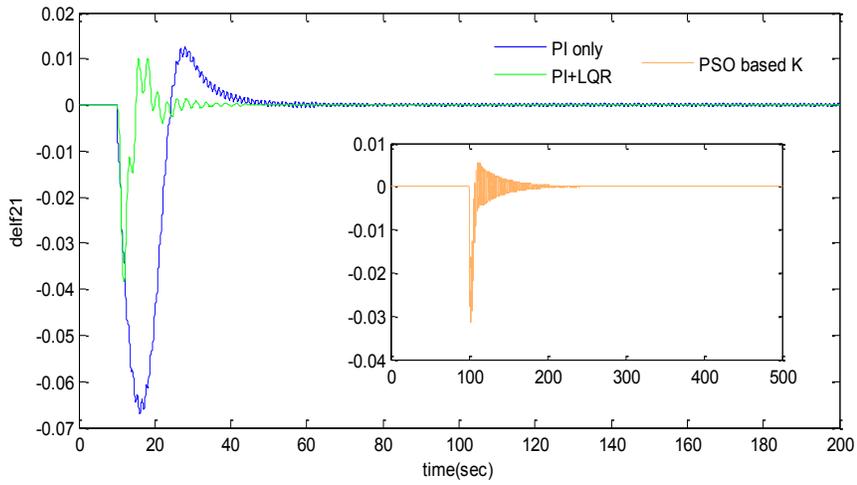

Fig. 4: Δf2 Vs t of area 2 if 1% load variation occurs in area 1

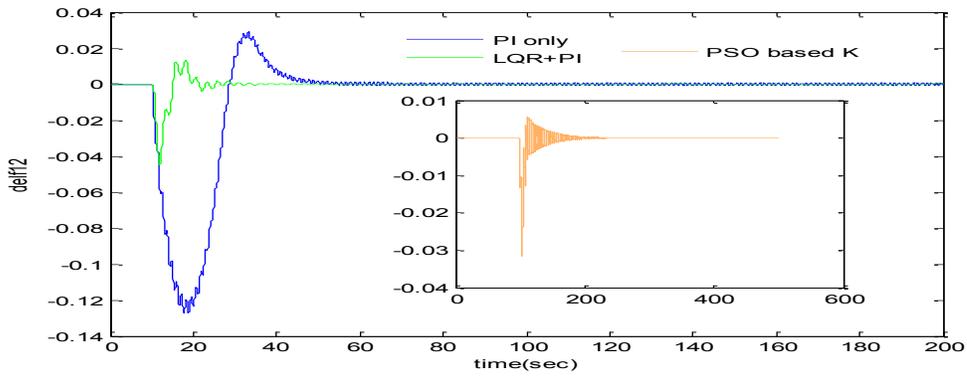

Fig. 5: Δf1 Vs t of area1 if 1% load variation occurs in area 2

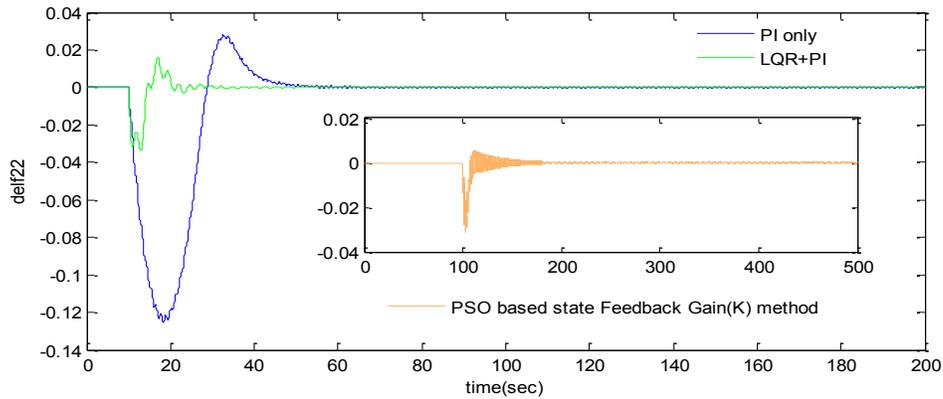

Fig. 6: Δf2 Vs t for area 2 if 1% load variation occurs in area 2

Based on the above results, the peak undershoot for the variation in frequency when there is 1% load change in one of areas in each case is shown in Table 1 and Table 2.The settling time for the responses is compared in Table3 and Table4.The peak overshoots are compared in Table 5 and table 6.The following nomenclature has been used for various frequency responses..

delf11: Frequency Response of power system Area-I for 1% load deviation in Area-I
delf21: Frequency Response of power system Area-II for 1% load deviation in Area-I
delf12: Frequency Response of power system Area-I for 1% load deviation in Area-II
delf22: Frequency Response of power system Area-II for 1% load deviation in Area-II

Table 1. Performance comparison of peak undershoot(Hz) with Integral , LQR-PI and PSO based state feedback (K) controller for power system network with eleven states in each area for 1% variation in load for Area-I

|  | Integral controller | LQR+PI controller | PSO based State feedback (K) controller |
|---|---|---|---|
| delf11 | -0.07 | -0.03 | -0.015 |
| delf21 | -0.07 | -0.04 | -0.03 |

Table 2. Performance comparison peak undershoot(Hz) with Integral, LQR-PI and PSO based state feedback (K) controller for power system network for 1% load variation in Area II

|  | Integral controller | LQR+PI controller | PSO based State feedback (K) controller |
|---|---|---|---|
| delf12 | -0.12 | -0.04 | -0.03 |
| delf22 | -0.12 | -0.03 | -0.015 |

Table 3. Performance comparison of settling time (seconds) with Integral, LQR-PI and PSO based state feedback (K) controller for power system network for 1% load variation in Area I

|  | Integral controller | LQR+PI controller | PSO based State feedback (K) controller |
|---|---|---|---|
| delf11 | 20s | 15s | 100s |
| delf21 | 18s | 15s | 100s |

Table 4. Performance comparison of settling time (seconds) with Integral, LQR-PI and PSO based state feedback (K) controller for power system network for 1% load variation in Area II

|  | Integral controller | LQR+PI controller | PSO based State feedback (K) controller |
|---|---|---|---|
| delf12 | 18s | 15s | 100s |
| delf22 | 18s | 15s | 100s |

Table 5. Performance comparison of peak overshoot (Hz) with Integral, LQR-PI and PSO based state feedback (K) controller for power system network for 1% load variation in Area I

|  | Integral controller | LQR+PI controller | PSO based State feedback (K) controller |
|---|---|---|---|
| delf11 | 0.011 | 0.011 | 0.005 |
| delf21 | 0.011 | 0.010 | 0.005 |

Table 6. Performance comparison of peak overshoot (Hz) with Integral, LQR-PI and PSO based state feedback (K) controller for power system network for 1% load variation in Area II

|        | Integral controller | LQR+PI controller | PSO based State feedback (K) controller |
|--------|---------------------|-------------------|-----------------------------------------|
| delf12 | 0.020               | 0.010             | 0.005                                   |
| delf22 | 0.022               | 0.015             | 0.005                                   |

The comparison of responses of three controllers for load frequency control shows that for same load variation the best performance is with PSO based state feedback (K) controller, then the LQR-PI controller followed by Integral controller performance for the peak overshoot and peak undershoot. As the number of variables are large in number and computational steps are more for optimization in PSO based state feedback (K) controller, so this method has more settling time followed by LQR-PI controller and Integral controller.

## 5. Conclusion and Future Work

In this work the performance comparison of three controllers such as PSO based state feedback gain (K) controller, LQR-PI controller and Integral controller performance has been done. The system investigated for frequency response enhancement is comprising of two thermal systems with reheat turbines. The two systems are interconnected with tie line. For the comparison of the performance of three types of controllers the transfer function model and state space model are developed for the system. Transfer function model has been used for finding the response of Integral controller and PSO based state feedback gain matrix (K) controller. Further, the state space based model is designed to observe the LQR-PI controller performance. The load variation in each area at one time is 1% and responses are obtained for each controller. The frequency variation is observed in both the areas although the load changes in one area only at one time. The frequency deviation is minimum when PSO based state feedback gain (K) controller is used followed by the LQR-PI controller. The maximum variation in frequency is in case of Integral controller which is although fully optimized with powerful computational technique PSO. This comparison of frequency controller is essentially useful in identifying the most suitable frequency controller which controls the frequency upon load changes for the interconnected power network. The present comparison gives the most efficient PSO based state feedback gain (K) controller which gives the minimum frequency deviation upon load variation in any area. In the present work the proposed technique of PSO based state feedback gain (K) controller for frequency control has been applied for two area thermal network, in future this novel technique can be applied for other types of interconnected power systems.